\documentclass[12pt]{amsart}
\usepackage{fullpage}
\usepackage[active]{srcltx}
\usepackage{amssymb,amsthm,amsmath,amscd}

\usepackage{mathtools}

\RequirePackage[dvipsnames,usenames]{xcolor}

\usepackage[colorlinks, linkcolor=Mahogany, citecolor=ForestGreen]{hyperref}

\input{kmacros3.sty}

\usepackage{tikz}
\usepackage{tikz-cd}
\usetikzlibrary{cd}
\usepackage{graphicx}
\usepackage[all,cmtip]{xy}

\usepackage{mabliautoref}
\numberwithin{equation}{theorem}

\usepackage{pifont}
\usepackage{upgreek}

\usepackage{ulem}
\usepackage{stmaryrd}

\numberwithin{equation}{theorem}

\usepackage{enumerate}
\usepackage{calc}

\usepackage{verbatim}
\usepackage{alltt}

\renewcommand{\m}{\mathfrak{m}}

\renewcommand{\fram}{\mathfrak{m}}

\newcommand{\perfd}{\textnormal{perfd}}

\newcommand{\Bl}{\textnormal{Bl}}

\newcommand{\Kos}{\textnormal{Kos}}
\newcommand{\trdeg}{\textnormal{trdeg}}

\begin{document}

\title{A short proof of Serre positivity}

\author{Linquan Ma}

\address{Department of Mathematics, Purdue University, West Lafayette, IN, USA}
\email{ma326@purdue.edu}

\begin{abstract}
Recently, OpenAI found a proof of the positivity of Serre's intersection multiplicity. The two main ingredients are a descent result on rational $K$-theory with support; and a construction of certain category of asymptotic complexes that identifies with the category of perfect complexes over a connective (non-commutative) derived ring. In this article, we present a simpler proof of Serre positivity that circumvents the second step; instead, we use a ``reduction to the diagonal" trick via the $A_{\text{inf}}$-construction. 
\end{abstract}

\maketitle

\section{Introduction}

In this paper, we explain a proof of the following theorem, which settles the long-standing Serre's positivity conjecture \cite{SerreBookLocalAlgebra1965} on intersection multiplicities in mixed characteristic. We fix a prime number $p$ throughout.

\begin{theorem}
\label{theorem: Serre positivity}
Let $(S,\m_S)$ be a ramified regular local ring of mixed characteristic $(0,p)$ and dimension $n$. Let $P$ and $Q$ be two prime ideals of $S$ such that $P+Q$ is $\m_S$-primary and $\dim(S/P) + \dim(S/Q) =n$. Then we have 
$$\chi(S/P, S/Q)>0.$$
\end{theorem}

We briefly recall some history and background on Serre's intersection multiplicity. Let $(S,\m_S)$ be a regular local ring and let $M$, $N$ be two finitely generated $S$-modules. Suppose $M\otimes_SN$ has finite length (i.e., $M\otimes_SN$ is only supported at $\{\m_S\}$). Serre introduced the following intersection multiplicity 
$$\chi(M, N) : = \chi(M\otimes^\myL_S N) = \sum_i (-1)^i\ell(\Tor_i^S(M,N)).$$

It was proved in \cite{SerreBookLocalAlgebra1965} that when the completion of $S$ is a formal power series ring over a field or a mixed characteristic DVR (this includes the equal characteristic and unramified mixed characteristic case), then 
\begin{enumerate}
    \item $\dim(M) + \dim(N) \leq \dim(S)$;
    \item $\chi(M, N)=0$ when $\dim(M) + \dim(N) <\dim(S)$;
    \item $\chi(M, N)>0$ when $\dim(M) + \dim(N) =\dim(S)$.
\end{enumerate}

Serre also proved that (a) holds for all regular local rings in \cite{SerreBookLocalAlgebra1965}, and conjectured    (``Il est naturel de ${\hbox{\rm conjecturer}}$")
that  (b) and (c) hold for all regular local rings as well. Roberts \cite{RobertsVanishing1985} and, independently, Gillet and Soul\'e \cite{GilletSouleVanishing1987}, proved part (b). Gabber proved that $\chi(M,N) \geq 0$, see \cite{BerthelotGabberNonnegativity,HochsterNonnegativity1997}. The remaining question is the strict positivity (i.e., part (c)) when $S$ is ramified mixed characteristic. By taking the prime cyclic filtrations of $M$ and $N$, and using the bi-additivity of $\chi(-,-)$ and the already established vanishing part (b), it is easy to see that part (c) reduces to \autoref{theorem: Serre positivity}. 

\subsection*{Outline of the proof} 
We first give an outline of the idea and strategy of the proof, which is largely inspired from \cite{BhattHochsterMaLimCM2026}. Our goal is to prove the positivity of $\chi(S/P\otimes^\myL_S S/Q)$. The first step is to replace $\chi$ by $\chi_\infty$, namely we show that 
$$\chi(S/P\otimes^\myL_S S/Q) = \chi_\infty^{\underline{x}}\big((S/P)^{A_\infty}_\perfd\otimes^\myL_S S/Q\big)$$
where $(S/P)^{A_\infty}_\perfd$ is a mixed characteristic analog of perfection in characteristic $p$ (see \cite[Section 8 and Section 10]{BhattScholzePrismaticCohomology2022}) and $\chi_\infty^{\underline{x}}$ is the alternating sum of the normalized length computed over $A_\infty$ (see \cite[Section 4]{CaiLeeMaSchwedeTuckerPerfdSignature2026}). This step is essentially the same as in \cite[Section 5.1, Lemma 5.2]{OpenAISerrePositivity2026} (both the strategy and the proof) and is carried out in \autoref{sec: first step}. But note that, instead of using the full results in \cite[Section 3]{OpenAISerrePositivity2026}, this step only requires a relatively simpler input from $K$-theory (\cite[Appendix A]{MaLech2026}), see the proof of \autoref{thm: chi equals chi_infty}. 

The second and more crucial step is to replace $\chi_\infty$ by $\chi_{\infty,\infty}$. More specifically, we construct $(S/Q)^{B_\infty}_\perfd$ similar to $(S/P)^{A_\infty}_\perfd$ and show that 
$$\chi_\infty^{\underline{x}}\big((S/P)^{A_\infty}_\perfd\otimes^\myL_S S/Q\big) =\chi_{\infty,\infty}^{\underline{x},\underline{y}}\big((S/P)^{A_\infty}_\perfd\otimes^\myL_S (S/Q)^{B_\infty}_\perfd\big).$$
Here $\chi_{\infty,\infty}^{\underline{x},\underline{y}}$ is the alternating sum of normalized length computed over some auxiliary perfectoid ring, see \autoref{setting: normalized length}. To establish this equality, we make use of an isogeny theorem of perfectoidization (see \cite[Proposition 2.4.7]{CaiLeeMaSchwedeTuckerPerfdSignature2026}) to reduce it to a vanishing statement of $\chi_{\infty,\infty}^{\underline{x},\underline{y}}((S/P)^{A_\infty}_\perfd\otimes^\myL_S G)$ where $G$ has smaller support, see \autoref{discuss: strategy of proof}. We then identify the intersection of $(S/P)^{A_\infty}_\perfd$ and $G$ over $S$ with the intersection of $A_{\inf}((S/P)^{A_\infty}_\perfd)$ and $G$ over an unramified regular local ring, which then allows us to apply reduction to the diagonal reinterpreting the intersection via a Koszul complex. The latter can be shown to vanish essentially because one of the diagonal elements is redundant, see \autoref{prop: analytic spread at most n-1} and the proof of \autoref{thm: vanishing and completion of proof}. 

Finally, by \cite[Corollary 4.0.4]{CaiLeeMaSchwedeTuckerPerfdSignature2026} (which is an immediate consequence of the almost purity theorem in the form of Bhatt--Scholze \cite[Theorem 10.9]{BhattScholzePrismaticCohomology2022}), we know that $(S/P)^{A_\infty}_\perfd$ and $(S/Q)^{B_\infty}_\perfd$ are almost Cohen--Macaulay. This implies that
$$\chi_{\infty,\infty}^{\underline{x},\underline{y}}\big((S/P)^{A_\infty}_\perfd\otimes^\myL_S (S/Q)^{B_\infty}_\perfd\big) = \lambda_{\infty,\infty}^{\underline{x},\underline{y}}\big((S/P)^{A_\infty}_\perfd\otimes_S (S/Q)^{B_\infty}_\perfd\big)$$
by \autoref{lem: strict positivity} (whose proof is similar in spirit to the argument in \cite[Theorem 6.1]{BhattHochsterMaLimCM2026}). The right hand side above is positive by \cite[Proposition 4.0.20]{CaiLeeMaSchwedeTuckerPerfdSignature2026}. %The second and third step diverge from OpenAI's proof.

\subsection*{Acknowledgement and AI disclosure} 
The overall strategy of this proof was originated during the preparation of \cite{BhattHochsterMaLimCM2026}. The author learned the results and methods in \autoref{sec: first step} from \cite[Section 3]{OpenAISerrePositivity2026} (and \cite[Appendix A]{MaLech2026}, which was discovered using ChatGPT). In addition, the author used ChatGPT to find the proof of \autoref{prop: analytic spread at most n-1}, as well as the arguments in \autoref{lem: extend valuation perfectoidization} and \autoref{lem: valuation on W(R)}.

I would like to thank Bhargav Bhatt for many enjoyable discussions on Serre's positivity conjecture and other applications of perfectoid techniques in commutative algebra. Many ideas and details of the proof presented here were discovered during our meetings and conversations with Bhargav. In addition, I want to thank Bhargav Bhatt and Mel Hochster for the collaboration \cite{BhattHochsterMaLimCM2026}. The author was supported by NSF Grant DMS \#2302430.

\newpage
\section{Positivity of Serre intersection multiplicity: the first step}
\label{sec: first step}

In this section we carry out the first step outlined in the introduction. We will work in the following setup. 

\begin{setting}
\label{setting: step 1}
Let $(S,\m_S,k)$ be a complete and ramified regular local ring of mixed characteristic $(0,p)$ of dimension $n$, with $k$ algebraically closed.  Let $P$, $Q$ be two prime ideals of $S$ such that $P+Q$ is $\m_S$-primary and $\dim(S/P)+\dim(S/Q)=n$. Suppose $S/P$ is $p$-torsion free. Set $R=S/P$ which is a Noetherian complete local domain.  Choose a system of parameters $\underline{x} :=p, x_2,\dots, x_d$ of $R$ and let $(A,\m_A):=W(k)[[x_2,\dots,x_d]]\subseteq R$ be the corresponding Noether--Cohen normalization. Set 
$$A_\infty:= A[p^{1/p^\infty}, x_2^{1/p^\infty},\dots,x_d^{1/p^\infty}]^{\wedge p}$$
and let $R^{A_\infty}_{\perfd}:= (R\otimes_AA_\infty)_{\perfd}$ where $(-)_{\perfd}$ is the perfectoidization functor of Bhatt--Scholze (see \cite[Section 8]{BhattScholzePrismaticCohomology2022}). For the convenience of future reference, we also set $A':=W(k)[x_2,\dots,x_d]^{\wedge p}$ and set
$$A'_\infty := W(k)[p^{1/p^\infty}, x_2^{1/p^\infty},\dots,x_d^{1/p^\infty}]^{\wedge p}.$$
\end{setting}

We will use the normalized length function $\lambda_\infty(-)$ of Faltings (see \cite{FaltingsAlmostEtale2002} or \cite[Section 2.2]{CaiLeeMaSchwedeTuckerPerfdSignature2026}), and we will sometimes use almost mathematics with respect to $(g)_\perfd$ (for some $0\neq g\in A$) as in \cite[Section 10]{BhattScholzePrismaticCohomology2022}. We point out that the essential properties we need are simply that there is a well-defined ``length function"
$$\lambda_\infty(-): \,\ \text{$\m_A$-power torsion $A_\infty$-modules} \longrightarrow \mathbb{R}_{\geq 0} \cup \{\infty\}$$
that is additive on short exact sequences, and that takes value zero on $A_\infty$-modules that are $g$-almost zero (i.e., annihilated by $(g)_\perfd$). Note that under the canonical map $A_\infty'\to A_\infty$, the normalized length computed over $A'_\infty$ and $A_\infty$ are the same: this is because the normalized length is defined on $\m_A$-power torsion modules. In the sequel, we will need to distinguish between the normalized length computed over various perfectoid rings, so from now on we will use $\lambda_\infty^{\underline{x}}(-)$ to denote the normalized length computed over $A_\infty$ or $A'_\infty$.

%We relate the Euler characteristic over $S$ with the ``normalized Euler characteristic" computed over $A_\infty$, when we replace $R=S/P$ by $R^{A_\infty}_\perfd$. 

\begin{theorem}
\label{thm: chi equals chi_infty}
With notations as in \autoref{setting: step 1}, for any finitely generated $S/Q$-module $N$ we have 
$$\chi(S/P\otimes^\myL_S N) = \chi_\infty^{\underline{x}}(R^{A_\infty}_\perfd \otimes^\myL_SN)$$
where the latter is defined as
$$\chi_\infty^{\underline{x}}(R^{A_\infty}_\perfd \otimes^\myL_SN) := \sum_i (-1)^i \lambda^{\underline{x}}_\infty\big(H_i(R^{A_\infty}_\perfd \otimes^\myL_SN)\big),$$
and each term in the summation is finite.
\end{theorem}
\begin{proof}
This follows from \cite[Section 2]{MaLech2026}. More specifically, let $F_\bullet$ be a finite free resolution of $N$ over $S$. Then we have $$\chi_\infty^{\underline{x}}(R^{A_\infty}_\perfd \otimes^\myL_SN) = \chi_\perfd^{\underline{x}}(F_\bullet \otimes^\myL_S R)$$ where the latter is defined as in \cite[Definition 2.4]{MaLech2026}. In particular, $\chi_\infty^{\underline{x}}$ is well-defined (i.e., each term in the summation is finite by \cite[Lemma 2.3]{MaLech2026}). Now by \cite[Theorem 2.9]{MaLech2026}, $$\chi_\perfd^{\underline{x}}(F_\bullet \otimes^\myL_S R) = \chi_\infty(F_\bullet \otimes^\myL_S R)$$ where $\chi_\infty$ denotes the Dutta multiplicity of the perfect complex $F_\bullet\otimes^\myL_SR$ over $R$ defined using localized Chern characters. Finally, it is well-known that under our setting,
$$\chi_\infty(F_\bullet \otimes^\myL_S R) = \chi(S/P\otimes^\myL_S N).$$
For example, this follows from \cite[calculation on the top of Page 376]{KuranoCharacteristicFreeDutta1993}.\footnote{Note that in \cite[Proposition 3.3]{KuranoCharacteristicFreeDutta1993}, the additional Cohen--Macaulay and characteristic hypotheses are only used for the positivity conclusion; those hypotheses are not needed for the Chern character calculation.} Putting these together completes the proof of the theorem.
\end{proof}

\begin{remark}
The proof of \autoref{thm: chi equals chi_infty} is essentially the same as in \cite{OpenAISerrePositivity2026}. The key behind the proof is a descent result on rational $K$-theory with support: see \cite[Theorem A.1]{MaLech2026} (which is used in \cite[Theorem 2.9]{MaLech2026} and is a special case of \cite[Theorem 3.1]{OpenAISerrePositivity2026}). This is the only place in our approach to Serre positivity that requires $K$-theory. Note that we have given an alternative and short exposition in \cite[Appendix A]{MaLech2026}.
\end{remark}

We next use \autoref{thm: chi equals chi_infty} and the existence of lim Cohen--Macaulay sequences in characteristic $p$ to prove Serre positivity when $S/Q$ has characteristic $p$.

\begin{theorem}
\label{thm: S/Q has char p}
With notations as in \autoref{setting: step 1}, assume further that $S/Q$ has characteristic $p$. Then $\chi(S/P \otimes^\myL_S S/Q)>0$.  
\end{theorem}
\begin{proof}
First note that there exists $K\gg 0$ so that 
$$\m_S^K\subseteq \ker\big(S\to \Hom_{D(S)}\big(R^{A_\infty}_\perfd\otimes^\myL_S F^e_*(S/Q), R^{A_\infty}_\perfd\otimes^\myL_S F^e_*(S/Q)\big)\big)$$
for all $e\geq 0$. This follows from a simple spectral sequence argument as $R^{A_\infty}_\perfd\otimes^\myL_S F^e_*(S/Q)\in D^{[-n,0]}(S)$ and each of its homology modules is annihilated by $P+Q$ (which is $\m_S$-primary and is independent of $e$). Fix a system of parameters $f_1,\dots,f_d, f_{d+1},\dots, f_n \in \m_S^K$ so that $f_1,\dots,f_d$ is a system of parameters on $R=S/P$ and $f_{d+1},\dots,f_n$ is a system of parameters on $S/Q$ (here we used the assumption that $\dim(S/P)+\dim(S/Q)=n$). The canonical map 
$$\Kos(f_1,\dots,f_n; R^{A_\infty}_\perfd\otimes^\myL_SF^e_*(S/Q)) \to \big(R^{A_\infty}_\perfd\otimes^\myL_SF^e_*(S/Q)\big)[n]$$
admits a section in $D(S)$. In particular, $\Tor_i^S(R^{A_\infty}_\perfd, F^e_*(S/Q))$ is a direct summand of 
\begin{align*}
& H_{i+n}\big(\Kos(f_1,\dots,f_n; R^{A_\infty}_\perfd \otimes^\myL_SF^e_*(S/Q))\big) \\
\cong & H_{i+n}\big(\Kos(f_1,\dots,f_d; R^{A_\infty}_\perfd) \otimes^\myL_S \Kos(f_{d+1},\dots, f_n; F^e_*(S/Q))\big).
\end{align*}
Examining the spectral sequence computing the latter, we find that its normalized length over $A_\infty$ is controlled by the sum of the normalized lengths of 
\begin{equation}
\label{eqn: Tor lim CM}
\Tor_{i+n-s-t}^S\big(H_s(\Kos(f_1,\dots,f_d; R^{A_\infty}_\perfd)), H_t\big(\Kos(f_{d+1},\dots, f_n; F^e_*(S/Q))\big)\big)
\end{equation}
for various $s$ and $t$. 

Note that if $s>0$, then the normalized length of \autoref{eqn: Tor lim CM} is zero by \cite[Proposition 4.0.10]{CaiLeeMaSchwedeTuckerPerfdSignature2026} since $H_s(\Kos(f_1,\dots,f_d; R^{A_\infty}_\perfd))$ is $g$-almost zero where $g$ is the discriminant of $A\to R$ by \cite[Corollary 4.0.4]{CaiLeeMaSchwedeTuckerPerfdSignature2026}; while if $s=0$ and $t>0$, then the normalized length of \autoref{eqn: Tor lim CM} is bounded by 
\begin{align*}
& \lambda_\infty^{\underline{x}}\big(\Tor_{i+n-t}^S(R^{A_\infty}_\perfd/(f_1,\dots,f_d)R^{A_\infty}_\perfd, k)\big)\cdot \ell\big(H_t\big(\Kos(f_{d+1},\dots, f_n; F^e_*(S/Q))\big)\big)  \\
\leq & \lambda_\infty^{\underline{x}}(R^{A_\infty}_\perfd/(f_1,\dots,f_d)R^{A_\infty}_\perfd)\cdot 2^n \cdot o(p^{e(n-d)})   
\end{align*}
since $F^e_*(S/Q)$ is a lim Cohen--Macaulay sequence by \cite[Theorem 5.4]{BhattHochsterMaLimCM2026}. Putting these together and noting that when $i>0$, \autoref{eqn: Tor lim CM} is nonzero only if $s>0$ or $t>0$, we have 
\begin{equation}
\label{eqn: eqn 1 in char p proof}
\lambda_\infty^{\underline{x}}\big(\Tor_i^S(R^{A_\infty}_\perfd, F^e_*(S/Q))\big)\leq \lambda_\infty^{\underline{x}}(R^{A_\infty}_\perfd/(f_1,\dots,f_d)R^{A_\infty}_\perfd)\cdot n2^n \cdot o(p^{e(n-d)}).
\end{equation}

On the other hand, when $i=0$, we have 
\begin{align}
\label{eqn: eqn 2 in char p proof}
\lambda_\infty^{\underline{x}}\big(R^{A_\infty}_\perfd \otimes_S F^e_*(S/Q)\big) & \geq \lambda_\infty^{\underline{x}}\big(R^{A_\infty}_\perfd/\m_SR^{A_\infty}_\perfd \otimes_k F^e_*(S/Q)/\m_SF^e_*(S/Q)\big) \\
& = \lambda_\infty^{\underline{x}}(R^{A_\infty}_\perfd/\m_SR^{A_\infty}_\perfd) \cdot \ell\big(F^e_*(S/Q)/\m_SF^e_*(S/Q)\big) \notag \\
& \geq \lambda_\infty^{\underline{x}}(R^{A_\infty}_\perfd/\m_SR^{A_\infty}_\perfd) \cdot \big(e_{\text{HK}}(S/Q)\cdot p^{e(n-d)} - o(p^{e(n-d)})\big) \notag \\
& \geq e_{\text{HK}}(S/Q)\cdot p^{e(n-d)} - o(p^{e(n-d)}), \notag
\end{align}
where the last inequality follows from \cite[Proposition 4.0.20]{CaiLeeMaSchwedeTuckerPerfdSignature2026}.

Finally, we have 
\begin{align*}
\chi(S/P \otimes^\myL_SS/Q) & = \lim_{e\to\infty} \frac{\chi\big(S/P \otimes^\myL_SF^e_*(S/Q)\big)}{p^{e(n-d)}} \\
& = \lim_{e\to\infty} \frac{\chi_\infty^{\underline{x}}\big(R^{A_\infty}_\perfd \otimes^\myL_SF^e_*(S/Q)\big)}{p^{e(n-d)}} \\
& \geq \lim_{e\to\infty} \frac{e_\text{HK}(S/Q)\cdot p^{e(n-d)} - o(p^{e(n-d)}) }{p^{e(n-d)}} \\
& \geq e_\text{HK}(S/Q)\geq 1,
\end{align*}
where the first equality follows from the vanishing part of Serre's intersection multiplicity (as $\rank(F^e_*(S/Q))=p^{e(n-d)}$), the second equality follows from \autoref{thm: chi equals chi_infty}; the first inequality follows from \autoref{eqn: eqn 1 in char p proof} and \autoref{eqn: eqn 2 in char p proof}. This completes the proof of the theorem.
\end{proof}

\section{Positivity of Serre intersection multiplicity: completion of the proof}
\label{sec: second step}

In this section we prove \autoref{theorem: Serre positivity} when both $S/P$ and $S/Q$ are $p$-torsion free (in fact, with small adjustments the argument also works when $S/Q$ has characteristic $p$, but we have already handled that case in \autoref{thm: S/Q has char p}). Our proof for this final step diverges significantly from the methods in \cite[Section 4 and Section 5]{OpenAISerrePositivity2026}. We start by setting up notations. 

\begin{setting}
\label{setting: normalized length}
We continue to use the notations as in \autoref{setting: step 1} and  assume further that $S/Q$ is $p$-torsion free and $\dim(S/Q)=c$. Choose a system of parameters $\underline{y}:= p, y_2,\dots,y_c$ of $S/Q$ and let $(B,\m_B):=W(k)[[y_2,\dots,y_c]]\subseteq S/Q$ be the corresponding Noether--Cohen normalization. Set $B':=W(k)[y_2,\dots,y_c]^{\wedge p}$ and set
$$B_\infty:= B[p^{1/p^\infty}, y_2^{1/p^\infty},\dots,y_c^{1/p^\infty}]^{\wedge p}, \,\ \,\  B'_\infty:= W(k)[p^{1/p^\infty}, y_2^{1/p^\infty},\dots,y_c^{1/p^\infty}]^{\wedge p}.$$ 
Let $(S/Q)^{B_\infty}_\perfd = (S/Q \otimes_BB_\infty)_\perfd$. 
The normalized length computed over $B_\infty$ and $B'_\infty$ are the same and will be denoted by $\lambda_\infty^{\underline{y}}(-)$. 
We will also compute normalized length over the following perfectoid ring:
$$C_\infty := W(k)[x_1^{1/p^\infty},x_2^{1/p^\infty}, \dots,x_d^{1/p^\infty}, p^{1/p^\infty}, y_2^{1/p^\infty},\dots,y_c^{1/p^\infty}]^{\wedge p}$$
and we denote the normalized length computed over $C_\infty$ by $\lambda_{\infty,\infty}^{\underline{x},\underline{y}}(-)$. Note that
$$C_\infty/(p-x_1) \cong (A'_\infty \otimes_{W(k)} B'_\infty)^{\wedge p}.$$
In particular, for any $\m_A$-power torsion $A_\infty$-module $M$ and $\m_B$-power torsion $B_\infty$-module $N$, $\lambda_{\infty,\infty}^{\underline{x},\underline{y}}(M \otimes_{W(k)} N)$ is well-defined. 
\end{setting}

\begin{remark}
We caution the readers that, with notations as in \autoref{setting: normalized length}, $(A'_\infty \otimes_{W(k)} B'_\infty)^{\wedge p}$ itself is not a perfectoid ring, and there are two compatible systems of $p$-power roots of $p$ in $(A'_\infty \otimes_{W(k)} B'_\infty)^{\wedge p}\cong C_\infty/(p-x_1)$, namely $\{x_1^{1/p^e}\}_{e\in\mathbb
N}$ and $\{p^{1/p^e}\}_{e\in\mathbb N}$. 
\end{remark}

We record some basic results relating $\lambda_{\infty,\infty}^{\underline{x},\underline{y}}(-)$ with $\lambda_\infty^{\underline{x}}(-)$ and $\lambda_\infty^{\underline{y}}(-)$.

\begin{lemma}
\label{lem: chi_inftyinfty leq product of chi_infty}
With notations as in \autoref{setting: normalized length}, for any $\m_A$-power torsion $A_\infty$-module $M$ and $\m_B$-power torsion $B_\infty$-module $N$, we have
$$\lambda_{\infty,\infty}^{\underline{x},\underline{y}}(M \otimes_{W(k)} N) \leq \lambda_\infty^{\underline{x}}(M) \lambda_\infty^{\underline{y}}(N)$$
with equality if both $M$ and $N$ are annihilated by $p$. 
\end{lemma}
\begin{proof}
We will safely view $M$ and $N$ as $A_\infty'$ and $B_\infty'$-modules. First suppose $M$ and $N$ are finitely presented over $A_\infty'$ and $B_\infty'$ respectively. There exists $e\gg0$ so that $M=M_e \otimes_{A'_e}A'_\infty$ and $N=N_e\otimes_{B'_e}B'_\infty$, where $M_e$ and $N_e$ are finitely generated over 
$$A_e':=W(k)[p^{1/p^e},x_2^{1/p^e},\dots,x_d^{1/p^e}]^{\wedge p} \,\ \text{and} \,\ B_e':=W(k)[p^{1/p^e},y_2^{1/p^e},\dots,y_c^{1/p^e}]^{\wedge p}$$
respectively. Then it is easy to see that 
$$M\otimes_{W(k)}N\cong (M_e\otimes_{W(k)}N_e)\otimes_{A'_e\otimes_{W(k)}B'_e}(A'_\infty\otimes_{W(k)}B'_\infty) \cong (M_e\otimes_{W(k)}N_e)\otimes_{C_e}C_\infty$$
where $C_e=W(k)[x_1^{1/p^e},\dots,x_d^{1/p^e}, p^{1/p^e},y_2^{1/p^e},\dots,y_c^{1/p^e}]^{\wedge p}$ and $C_e/(p-x_1)\cong(A'_e\otimes_{W(k)}B'_e)^{\wedge p}$. It follows from the definition of normalized length that 
$$\lambda_\infty^{\underline{x}}(M)= \frac{1}{p^{ed}}\ell(M_e), \,\, \lambda_\infty^{\underline{y}}(N)= \frac{1}{p^{ec}}\ell(N_e), \,\ \lambda_{\infty,\infty}^{\underline{x},\underline{y}}(M\otimes_{W(k)}N)= \frac{1}{p^{en}}\ell(M_e\otimes_{W(k)}N_e).$$
Thus the conclusion follows since $d+c=n$ and $\ell(M_e\otimes_{W(k)}N_e) \leq \ell(M_e)\ell(N_e)$, 
with equality when $M_e$ and $N_e$ are annihilated by $p$ (which holds when $M$ and $N$ are annihilated by $p$).

Next we assume $M$ and $N$ are finitely generated over $A_\infty'$ and $B_\infty'$ respectively. We write $M=\colim_i M_i$ and $N=\colim_i N_i$ so that $M_i$ and $N_i$ are finitely presented over $A_\infty'$ and $B_\infty'$ respectively; and the transition maps are surjective. It follows that 
\begin{align*}
\lambda_{\infty,\infty}^{\underline{x},\underline{y}}(M\otimes_{W(k)}N) = \lim_i \lambda_{\infty,\infty}^{\underline{x},\underline{y}}(M_i \otimes_{W(k)} N_i) \leq \lim_i \lambda_\infty^{\underline{x}}(M_i)\lambda_\infty^{\underline{y}}(N_i) = \lambda_\infty^{\underline{x}}(M)\lambda_\infty^{\underline{y}}(N)
\end{align*}
where the first and last equalities follow from \cite[Proposition 2.2.6 (b)]{CaiLeeMaSchwedeTuckerPerfdSignature2026} and the middle inequality follows from the already established case when $M_i$ and $N_i$ are finitely presented; and this inequality is an equality when $M_i$ and $N_i$ are annihilated by $p$ (which we may arrange when $M$ and $N$ are annihilated by $p$). 

Finally, for general $M$ and $N$, we write $M=\colim_i M_i$ and $N=\colim_i N_i$ so that $M_i$ and $N_i$ are finitely generated over $A_\infty'$ and $B_\infty'$ respectively; and the transition maps are injective. Consider the following diagram for all $i\leq j$:
\[\xymatrix{
M_i \otimes_{W(k)} N_i \ar@{->>}[r] \ar[d] & \Image(i) \ar@{^{(}->}[r] \ar@{^{(}->}[d] & M \otimes_{W(k)} N \ar@{=}[d]\\
M_j \otimes_{W(k)} N_j \ar@{->>}[r]  & \Image(j) \ar@{^{(}->}[r]  & M \otimes_{W(k)} N. \\
}\]
It follows from this diagram that 
\begin{align*}
\lambda_{\infty,\infty}^{\underline{x},\underline{y}}(M\otimes_{W(k)}N) & = \lim_i \lambda_{\infty,\infty}^{\underline{x},\underline{y}}(\Image(i)) \\
& \leq  \lim_i \lambda_{\infty,\infty}^{\underline{x},\underline{y}}(M_i\otimes_{W(k)}N_i)  \\
& \leq \lim_i \lambda_\infty^{\underline{x}}(M_i)\lambda_\infty^{\underline{y}}(N_i) \\
& = \lambda_\infty^{\underline{x}}(M)\lambda_\infty^{\underline{y}}(N).
\end{align*}
Here the first and last equalities follow from \cite[Proposition 2.2.6 (b)]{CaiLeeMaSchwedeTuckerPerfdSignature2026} (as the transition maps are injective); the first inequality follows from the surjection $M_i \otimes_{W(k)} N_i \twoheadrightarrow \Image(i)$ and the second inequality follows from the already established case when $M_i$ and $N_i$ are finitely generated. Note that when $M$ and $N$ are annihilated by $p$, the first inequality is an equality: both $M_i$ and $N_i$ are annihilated by $p$ and thus $M_i \otimes_{W(k)} N_i \cong M_i \otimes_{k} N_i$ and $M\otimes_{W(k)} N \cong M \otimes_{k} N$; which implies that $M_i \otimes_{W(k)} N_i \cong \Image(i)$. Since the second inequality is also an equality by the already established case when $M_i$ and $N_i$ are finitely generated and annihilated by $p$, we have equalities throughout. 
\end{proof}

\begin{lemma}
\label{lem: chi_inftyinfty = chi_infty trivial extension}
With notations as in \autoref{setting: normalized length}. Suppose $M$ is a module over $$W(k)[x_1^{1/p^\infty},x_2^{1/p^\infty}, \dots,x_d^{1/p^\infty},y_2,\dots,y_c]^{\wedge p}/(p-x_1)$$ that is $(x_1,\dots,x_d, p, y_2,\dots,y_c)$-power torsion. Then we have 
$$\lambda_{\infty,\infty}^{\underline{x},\underline{y}}(M\otimes_{B'}B'_\infty)=\lambda_\infty^{\underline{x}}(M).$$
\end{lemma}
\begin{proof}
By writing $M=\colim_i M_i$ with injective transition maps and with $M_i$ finitely generated over $W(k)[x_1^{1/p^\infty},x_2^{1/p^\infty}, \dots,x_d^{1/p^\infty},y_2,\dots,y_c]^{\wedge p}/(p-x_1)$, it follows immediately from \cite[Proposition 2.2.6 (b)]{CaiLeeMaSchwedeTuckerPerfdSignature2026} that it is enough to prove the lemma for each $M_i$ (note that $B'_\infty$ is $p$-completely faithfully flat over $B'$ and thus $M\otimes_{B'}B'_\infty \cong \colim_i (M_i\otimes_{B'}B'_\infty)$ with injective transition maps). Thus we may assume $M$ is finitely generated. But then by writing $M=\colim_i M_i$ with surjective transition maps and with $M_i$ finitely presented over $W(k)[x_1^{1/p^\infty},x_2^{1/p^\infty}, \dots,x_d^{1/p^\infty},y_2,\dots,y_c]^{\wedge p}/(p-x_1)$, another use of \cite[Proposition 2.2.6 (b)]{CaiLeeMaSchwedeTuckerPerfdSignature2026}  reduces the lemma to the case that $M$ is finitely presented. We now write 
\begin{align*}
 M & = M_e \otimes_{W(k)[x_1^{1/p^e}, \dots,x_d^{1/p^e},y_2,\dots,y_c]^{\wedge p}/(p-x_1)}W(k)[x_1^{1/p^\infty}, \dots,x_d^{1/p^\infty},y_2,\dots,y_c]^{\wedge p}/(p-x_1) \\   
 & = M_e \otimes_{W(k)[x_1^{1/p^e}, \dots,x_d^{1/p^e}]^{\wedge p}/(p-x_1)}W(k)[x_1^{1/p^\infty},\dots,x_d^{1/p^\infty}]^{\wedge p}/(p-x_1) \\
 & \cong  M_e \otimes_{A_e'} A_\infty'.
\end{align*}
Thus by definition of normalized length over $A_\infty'$, we have 
$$\lambda_\infty^{\underline{x}}(M) = \frac{1}{p^{ed}}\ell(M_e).$$

On the other hand, we have 
\begin{align*}
M \otimes_{B'}B_\infty' & = M_e \otimes_{W(k)[x_1^{1/p^e}, \dots,x_d^{1/p^e},y_2,\dots,y_c]^{\wedge p}}W(k)[x_1^{1/p^\infty}, \dots,x_d^{1/p^\infty},y_2,\dots,y_c]^{\wedge p} \otimes_{B'}B_\infty' \\
& = M_e \otimes_{W(k)[x_1^{1/p^e}, \dots,x_d^{1/p^e},y_2,\dots,y_c]^{\wedge p}}W(k)[x_1^{1/p^\infty}, \dots,x_d^{1/p^\infty},p^{1/p^\infty}, y_2^{1/p^\infty},\dots,y_c^{1/p^\infty}]^{\wedge p}\\
& =  M_e \otimes_{W(k)[x_1^{1/p^e}, \dots,x_d^{1/p^e},y_2,\dots,y_c]^{\wedge p}} C_\infty \\
& = \big(M_e \otimes_{W(k)[x_1^{1/p^e}, \dots,x_d^{1/p^e},y_2,\dots,y_c]^{\wedge p}} C_e\big)\otimes_{C_e}C_\infty.
\end{align*}
It follows from the definition of normalized length over $C_\infty$ that 
\begin{align*}
    \lambda_{\infty,\infty}^{\underline{x},\underline{y}}(M \otimes_{B'}B_\infty') & = \frac{1}{p^{en}} \ell\big(M_e \otimes_{W(k)[x_1^{1/p^e}, \dots,x_d^{1/p^e},y_2,\dots,y_c]^{\wedge p}} C_e\big) \\
    & = \frac{1}{p^{en}} \cdot p^{ec} \cdot \ell(M_e) \\
    & = \frac{1}{p^{ed}} \ell(M_e)
\end{align*}
where the second equality follows from that $C_e= W(k)[x_1^{1/p^e},\dots,x_d^{1/p^e},p^{1/p^e}, y_2^{1/p^e},\dots,y_c^{1/p^e}]^{\wedge p}$ is $p$-completely free of rank $p^{ec}$ over $W(k)[x_1^{1/p^e}, \dots,x_d^{1/p^e},y_2,\dots,y_c]^{\wedge p}$ and the third equality follows from that $d+c=n$. Putting the above together completes the proof of the lemma.
\end{proof}

The next simple but crucial lemma is essentially the reason we could apply the reduction to the diagonal trick later.

\begin{lemma}
\label{lem: Tor-independent}
With notations as in \autoref{setting: step 1} and \autoref{setting: normalized length} and write $S=T/(p-f)$ where $T= W(k)[[z_1,\dots,z_n]]$ with $f\in (p,z_1,\dots,z_n)^2$. There is a map $T\to A_{\inf}(R^{A_\infty}_\perfd)$ so that the following diagram is a Tor-independent pushout:
\[
\xymatrix{
T \ar[d] \ar[r] & A_{\inf}(R^{A_\infty}_\perfd) \ar[d] \\
S \ar[r] & R^{A_\infty}_\perfd.
}
\]
Fixing such a map $T\to A_{\inf}(R^{A_\infty}_\perfd)$, we have a commutative diagram of natural maps:
\[\xymatrix@C-=.5pc{
C_\infty \ar[r] \ar[d] & (A_{\inf}(A_\infty') \otimes_{W(k)}B_\infty')^{\wedge p} \ar[r] \ar[d] & (A_{\inf}(R^{A_\infty}_\perfd) \otimes_{W(k)} B_\infty)^{\wedge p}  \ar[d]   \\
C_\infty/(p-x_1) \ar[r]^-{\sim} & (A_\infty' \otimes_{W(k)} B_\infty')^{\wedge p} \ar[d] & A_{\inf}(R^{A_\infty}_\perfd) \otimes_T(S/Q \otimes_BB_\infty) \ar@{=}[d]  \\
& (A_\infty \otimes_{W(k)} B_\infty)^{\wedge p} \ar[r] & R^{A_\infty}_\perfd \otimes_S(S/Q\otimes_BB_\infty) \ar[r] & R^{A_\infty}_\perfd \otimes_S(S/Q)^{B_\infty}_\perfd.
}\]
In particular, for any module $\mathcal{M}$ over $R^{A_\infty}_\perfd \otimes_S(S/Q)^{B_\infty}_\perfd$, the normalized length $\lambda_{\infty, \infty}^{\underline{x}, \underline{y}}(\mathcal{M})$ is well-defined.
\end{lemma}
\begin{proof}
Choose elements $\widetilde{z}_1,\dots,\widetilde{z}_n$ in $A_{\inf}(R^{A_\infty}_\perfd)$ that are arbitrary lifts of $z_1,\dots,z_n$ in $R^{A_\infty}_\perfd$. We obtain a map 
$$W(k)[z_1,\dots,z_n]^{\wedge p} \to A_{\inf}(R^{A_\infty}_\perfd)$$
sending $z_i$ to $\widetilde{z}_i$,
so that after modulo $p^m$ it fits into the following commutative diagram:
\[\xymatrix{
W(k)[z_1,\dots,z_n]/p^m \ar[r] \ar[d] & A_{\inf}(R^{A_\infty}_\perfd)/p^m \ar[r] & A_{\inf}(R^{A_\infty}_\perfd)/(p,\xi)^m \ar[d] \\
T/p^m\cong W(k)[[z_1,\dots,z_n]]/p^m \ar[r] \ar@{.>}[rru]^{\alpha_m} & R^{A_\infty}_\perfd/p^m \cong A_{\inf}(R^{A_\infty}_\perfd)/(p^m,\xi) \ar[r] & A_{\inf}(R^{A_\infty}_\perfd)/(p,\xi) 
}\]
where $(\xi)=\ker(A_{\inf}(R^{A_\infty}_\perfd) \to R^{A_\infty}_\perfd)$. We claim that there is a unique map $\alpha_m$ (exhibited by the dashed arrow) making the diagram above commute. To see this claim, it is enough to observe that
$$W(k)/p^m[z_1,\dots,z_n] \to W(k)/p^m[[z_1,\dots,z_n]]$$
is formally \'etale for all $m$: for $m=1$ the relative Frobenius of $k[z_1,\dots,z_n] \to k[[z_1,\dots,z_n]]$ is an isomorphism (we are assuming $k$ is algebraically closed) so the relative cotangent complex vanishes, which then implies that for all $m$, the relative cotangent complex after derived mod $p$ vanishes so the relative cotangent complex itself vanishes by derived Nakayama. Thus we can take inverse limit of $\alpha_m$ to obtain a map
$$T=\varprojlim_m T/p^m \to \varprojlim_m A_{\inf}(R^{A_\infty}_\perfd)/(p,\xi)^m = A_{\inf}(R^{A_\infty}_\perfd).$$

To show the Tor-independency of the diagram, by \cite[Lemma 2.24]{BhattScholzePrismaticCohomology2022} it suffices to show that $p-f$ maps to a distinguished element in $A_{\inf}(R^{A_\infty}_\perfd)$. But as $f\in (p,z_1,\dots,z_n)^2$, it is straightforward to check that $\delta(p-f)\equiv\delta(p)$ mod $(p,z_1,\dots,z_n)A_{\inf}(R^{A_\infty}_\perfd)$ and thus $\delta(p-f)$ is a unit in $A_{\inf}(R^{A_\infty}_\perfd)$.

Finally, the map $C_\infty \to (A_{\inf}(A_\infty') \otimes_{W(k)}B_\infty')^{\wedge p}$ is constructed by sending $x_1\mapsto [p^\flat]\in A_{\inf}(A_\infty')$ and $x_i\mapsto [x_i^\flat] \in A_{\inf}(A_\infty')$ for $i\geq 2$. It is straightforward to check that after we modulo $(p-x_1)$ this map agrees with $C_\infty/(p-x_1) \to (A_\infty'\otimes_{W(k)}B_\infty')^{\wedge p}$. The rest of the commutative diagram is essentially obvious. For the last assertion, note that any module $\mathcal{M}$ over $R^{A_\infty}_\perfd \otimes_S(S/Q)^{B_\infty}_\perfd$ is $\m_S$-power torsion, and thus $(x_1,\dots,x_d,p,y_2,\dots,y_c)$-power torsion when viewed as a $C_\infty$-module via the commutative diagram.
\end{proof}

The next lemma is basically saying that $\chi_{\infty, \infty}^{\underline{x}, \underline{y}}(- \otimes_S^\myL -)$ of two almost Cohen--Macaulay modules is positive. This should be viewed as the almost analog of the classical fact that $\chi(- \otimes^\myL_S -)$ of two Cohen--Macaulay modules (whose dimensions add up correctly) is positive. The proof is largely inspired from the theory of lim Cohen--Macaulay sequences developed in \cite{BhattHochsterMaLimCM2026} (once we set up $\lambda_{\infty, \infty}^{\underline{x}, \underline{y}}(-)$ and $\chi_{\infty, \infty}^{\underline{x}, \underline{y}}(-)$ correctly).

\begin{lemma}
\label{lem: strict positivity}
With notations as in \autoref{setting: step 1} and \autoref{setting: normalized length}. Then we have 
$$
\begin{cases}
\lambda_{\infty, \infty}^{\underline{x}, \underline{y}}\big(\Tor_i^S(R^{A_\infty}_\perfd, (S/Q)^{B_\infty}_\perfd)\big)  =0 \text{ when $i>0$;}\\
\lambda_{\infty, \infty}^{\underline{x}, \underline{y}}\big(\Tor_i^S(R^{A_\infty}_\perfd, (S/Q)^{B_\infty}_\perfd)\big)  >0 \text{ when $i=0$.}
\end{cases}
$$
In particular,
$$\chi_{\infty, \infty}^{\underline{x}, \underline{y}}\big(R^{A_\infty}_\perfd \otimes^\myL_S(S/Q)^{B_\infty}_\perfd\big) = \sum_i(-1)^i\lambda_{\infty, \infty}^{\underline{x}, \underline{y}}\big(\Tor_i^S(R^{A_\infty}_\perfd, (S/Q)^{B_\infty}_\perfd)\big) >0.$$
\end{lemma}
\begin{proof}
We first note that there exists $K\gg0$ so that
$$\m_S^K\subseteq \ker\big(S \to \Hom_{D(S)}\big(R^{A_\infty}_\perfd \otimes^\myL_S(S/Q)^{B_\infty}_\perfd, R^{A_\infty}_\perfd \otimes^\myL_S(S/Q)^{B_\infty}_\perfd\big)\big).$$
This follows from a simple spectral sequence argument as $R^{A_\infty}_\perfd \otimes^\myL_S(S/Q)^{B_\infty}_\perfd\in D^b(S)$ and each of its homology modules is annihilated by $P+Q$ which is $\m_S$-primary. We now fix a system of parameters of $f_1,\dots,f_d, f_{d+1},\dots, f_n \in\m_S^K \subseteq P+Q$ of $S$ such that $f_1,\dots,f_d$ is a system of parameters on $R=S/P$ and $f_{d+1},\dots,f_{n}$ is a system of parameters on $S/Q$. The canonical map
$$\Kos(f_1,\dots,f_n; R^{A_\infty}_\perfd \otimes^\myL_S(S/Q)^{B_\infty}_\perfd) \to \big(R^{A_\infty}_\perfd \otimes^\myL_S(S/Q)^{B_\infty}_\perfd\big)[n]$$
admits a section in $D(S)$. It follows that 
$\Tor_i^S(R^{A_\infty}_\perfd, (S/Q)^{B_\infty}_\perfd)$ is a direct summand of 
\begin{align*}
\mathcal{M}_i := & H_{i+n}\big(\Kos(f_1,\dots,f_n; R^{A_\infty}_\perfd \otimes^\myL_S(S/Q)^{B_\infty}_\perfd)\big) \\
\cong & H_{i+n}\big(\Kos(f_1,\dots,f_d; R^{A_\infty}_\perfd) \otimes^\myL_S \Kos(f_{d+1},\dots, f_n; (S/Q)^{B_\infty}_\perfd)\big).
\end{align*}
Examining the spectral sequence computing $\mathcal{M}_i$, we know that $\lambda_{\infty, \infty}^{\underline{x}, \underline{y}}(\mathcal{M}_i)$ is controlled by the sum of $\lambda_{\infty, \infty}^{\underline{x}, \underline{y}}(\mathcal{M}_{i,s,t})$, where 
$$\mathcal{M}_{i,s,t} :=\Tor_{i+n-s-t}^S\big(H_s(\Kos(f_1,\dots,f_d; R^{A_\infty}_\perfd)), H_t(\Kos(f_{d+1},\dots, f_n; (S/Q)^{B_\infty}_\perfd))\big).$$

Suppose $i>0$, then $\mathcal{M}_{i,s,t}\neq 0$ only if $s\geq 1$ or $t\geq 1$. Thus either $H_s\big(\Kos(f_1,\dots,f_d; R^{A_\infty}_\perfd)\big)$ is $g$-almost zero or $H_t\big(\Kos(f_{d+1},\dots, f_n; (S/Q)^{B_\infty}_\perfd)\big)$ is $h$-almost zero, where $g$ and $h$ are the discriminants of $A\to R$ and $B\to S/Q$ respectively by \cite[Corollary 4.0.4]{CaiLeeMaSchwedeTuckerPerfdSignature2026}. Therefore, if we view $\mathcal{M}_{i,s,t}$ as a module over $A_\infty\otimes_{W(k)}B_\infty$ via 
$$A_\infty\otimes_{W(k)}B_\infty \to R^{A_\infty}_\perfd\otimes_S(S/Q)^{B_\infty}_\perfd,$$
then either $\mathcal{M}_{i,s,t}$ is $g$-almost zero as an $A_\infty$-module or $\mathcal{M}_{i,s,t}$ is $h$-almost zero as a $B_\infty$-module. Since $\mathcal{M}_{i,s,t}$ is annihilated by $P+Q$ which is $\m_S$-primary, there exists $L\gg0$ so that $\m_A^L\mathcal{M}_{i,s,t}=0=\m_B^L\mathcal{M}_{i,s,t}$. Now if $\mathcal{M}_{i,s,t}$ is $g$-almost zero as an $A_\infty$-module, then by the $A_\infty\otimes_{W(k)}B_\infty$-linear surjection
$$\mathcal{M}_{i,s,t} \otimes_{W(k)}B_\infty/\m_B^LB_\infty \twoheadrightarrow\mathcal{M}_{i,s,t},$$
we have that
$$\lambda_{\infty, \infty}^{\underline{x}, \underline{y}}(\mathcal{M}_{i,s,t}) \leq \lambda_{\infty, \infty}^{\underline{x}, \underline{y}}(\mathcal{M}_{i,s,t} \otimes_{W(k)}B_\infty/\m_B^LB_\infty)\leq \lambda_\infty^{\underline{x}}(\mathcal{M}_{i,s,t})\cdot \lambda_\infty^{\underline{y}}(B_\infty/\m_B^LB_\infty)=0$$
where the second inequality follows from \autoref{lem: chi_inftyinfty leq product of chi_infty} and the equality follows from \cite[Proposition 4.0.10]{CaiLeeMaSchwedeTuckerPerfdSignature2026}. Similarly $\lambda_{\infty, \infty}^{\underline{x}, \underline{y}}(\mathcal{M}_{i,s,t})=0$ if $\mathcal{M}_{i,s,t}$ is $h$-almost zero as a $B_\infty$-module. Putting these together, we have 
$$\lambda_{\infty, \infty}^{\underline{x}, \underline{y}}\big(\Tor_i^S(R^{A_\infty}_\perfd, (S/Q)^{B_\infty}_\perfd)\big)\leq \lambda_{\infty, \infty}^{\underline{x}, \underline{y}}(\mathcal{M}_{i}) \leq \sum_{s,t} \lambda_{\infty, \infty}^{\underline{x}, \underline{y}}(\mathcal{M}_{i,s,t}) =0$$
for all $i>0$. Finally, when $i=0$, we have natural surjections
$$\big(R^{A_\infty}_\perfd/\m_S^K R^{A_\infty}_\perfd\big)\otimes_{W(k)}\big((S/Q)^{B_\infty}_\perfd/\m^K_S(S/Q)^{B_\infty}_\perfd\big) \twoheadrightarrow R^{A_\infty}_\perfd\otimes_S(S/Q)^{B_\infty}_\perfd$$
$$R^{A_\infty}_\perfd\otimes_S(S/Q)^{B_\infty}_\perfd \twoheadrightarrow \big(R^{A_\infty}_\perfd/\m_S R^{A_\infty}_\perfd\big)\otimes_{k}\big((S/Q)^{B_\infty}_\perfd/\m_S(S/Q)^{B_\infty}_\perfd\big).$$
The first surjection, \autoref{lem: chi_inftyinfty leq product of chi_infty}, and \cite[Corollary 4.0.18]{CaiLeeMaSchwedeTuckerPerfdSignature2026} imply that $\lambda_{\infty, \infty}^{\underline{x}, \underline{y}}\big(R^{A_\infty}_\perfd\otimes_S(S/Q)^{B_\infty}_\perfd\big)<\infty$; while the second surjection, \autoref{lem: chi_inftyinfty leq product of chi_infty}, and \cite[Proposition 4.0.20]{CaiLeeMaSchwedeTuckerPerfdSignature2026} imply that $\lambda_{\infty, \infty}^{\underline{x}, \underline{y}}\big(R^{A_\infty}_\perfd\otimes_S(S/Q)^{B_\infty}_\perfd\big)>0$.
%$$\lambda_{\infty, \infty}^{\underline{x}, \underline{y}}\big(R^{A_\infty}_\perfd\otimes_S(S/Q)^{B_\infty}_\perfd\big) \geq \lambda_\infty^{\underline{x}}(R^{A_\infty}_\perfd/\m_S R^{A_\infty}_\perfd) \cdot \lambda_\infty^{\underline{y}}\big((S/Q)^{B_\infty}_\perfd/\m_S(S/Q)^{B_\infty}_\perfd\big)\geq 1$$
%where the last inequality follows from \cite[Proposition 4.0.20]{CaiLeeMaSchwedeTuckerPerfdSignature2026}. This completes the proof of the lemma.
\end{proof}

%The next lemma relates $\chi_\infty^{\underline{x}}(-)$ with $\chi_{\infty,\infty}^{\underline{x},\underline{y}}(-)$.

\begin{lemma}
\label{lem: trivial extension}
With notations as in \autoref{setting: step 1} and \autoref{setting: normalized length}. Then we have 
$$\chi_{\infty,\infty}^{\underline{x},\underline{y}}\big(R^{A_\infty}_\perfd \otimes^\myL_S((S/Q)\otimes_BB_\infty)\big)= \chi_\infty^{\underline{x}}\big(R^{A_\infty}_\perfd \otimes^\myL_S(S/Q)\big).$$
\end{lemma}
\begin{proof}
Since $B_\infty$ is faithfully flat over $B$ and $R^{A_\infty}_\perfd \otimes^\myL_S(S/Q)$ is $\m_S$-power torsion, we have 
\begin{equation}
\label{eqn: eqn 1 in proof of trivial extension}
R^{A_\infty}_\perfd \otimes^\myL_S((S/Q)\otimes_BB_\infty) = (R^{A_\infty}_\perfd \otimes^\myL_S(S/Q))\otimes^\myL_BB_\infty = (R^{A_\infty}_\perfd \otimes^\myL_S(S/Q))\otimes^\myL_{B'}B'_\infty.
\end{equation}
Now each homology of $R^{A_\infty}_\perfd \otimes^\myL_S(S/Q)$ is a module over $R^{A_\infty}_\perfd \otimes_S(S/Q)$ that is $\m_S$-power torsion, and in particular a module over 
$$W(k)[x_1^{1/p^\infty},x_2^{1/p^\infty}, \dots,x_d^{1/p^\infty},y_2,\dots,y_c]^{\wedge p}/(p-x_1)$$ that is $(x_1,\dots,x_d, p, y_2,\dots,y_c)$-power torsion (see the last commutative diagram in \autoref{lem: Tor-independent}). By \autoref{lem: chi_inftyinfty = chi_infty trivial extension}, we have that 
\begin{align}
\label{eqn: eqn 2 in proof of trivial extension}
\lambda_\infty^{\underline{x}}\big(H_i(R^{A_\infty}_\perfd \otimes^\myL_S(S/Q))\big) & = \lambda_{\infty,\infty}^{\underline{x},\underline{y}}\big(H_i(R^{A_\infty}_\perfd \otimes^\myL_S(S/Q))\otimes_{B'}B'_\infty\big) \\
& = \lambda_{\infty,\infty}^{\underline{x},\underline{y}}\big(H_i((R^{A_\infty}_\perfd \otimes^\myL_S(S/Q))\otimes^\myL_{B'}B'_\infty)\big)    \notag
\end{align}
where the second equality follows from the $p$-complete flatness of $B_\infty'$ over $B'$ (in particular this shows that the right hand side of \autoref{eqn: eqn 2 in proof of trivial extension} is finite, as this is so for the left hand side by \autoref{thm: chi equals chi_infty}). Now taking the alternating sum on both sides of \autoref{eqn: eqn 2 in proof of trivial extension} and using \autoref{eqn: eqn 1 in proof of trivial extension} completes the proof of the lemma.
\end{proof}

We collect what we have proved so far, outline the strategy, and point out the remaining last step.

\begin{discussion}
\label{discuss: strategy of proof}
With notations as in \autoref{setting: step 1} and \autoref{setting: normalized length}, by the isogeny theorem of perfectoidization \cite[Proposition 2.4.7]{CaiLeeMaSchwedeTuckerPerfdSignature2026}, we have a short exact sequence:
\begin{equation}
\label{eqn: isogeny}
0 \to (S/Q)\otimes_BB_\infty \to (S/Q)^{B_\infty}_\perfd \to G \to 0
\end{equation}
where $G$ is annihilated by $h^K$ for some $K>0$ and $h\notin Q$ is the discriminant of $B\to (S/Q)$. Replacing $h$ by $h^K$, we have that $G$ is a module over $S/(Q+h) \otimes_BB_\infty$. By \autoref{thm: chi equals chi_infty} and \autoref{lem: trivial extension}, we have 
\begin{align*}
\chi(S/P\otimes^\myL_SS/Q) & = \chi_\infty^{\underline{x}}\big(R^{A_\infty}_\perfd \otimes^\myL_S(S/Q)\big) \\
& =\chi_{\infty,\infty}^{\underline{x},\underline{y}}\big(R^{A_\infty}_\perfd \otimes^\myL_S((S/Q)\otimes_BB_\infty)\big) \\
&  =  \chi_{\infty, \infty}^{\underline{x}, \underline{y}}\big(R^{A_\infty}_\perfd \otimes^\myL_S(S/Q)^{B_\infty}_\perfd\big) - \chi_{\infty, \infty}^{\underline{x}, \underline{y}}(R^{A_\infty}_\perfd \otimes^\myL_SG)
\end{align*}
where the last equality follows from \autoref{eqn: isogeny}. By \autoref{lem: strict positivity}, we will have $\chi(S/P\otimes^\myL_SS/Q)>0$, provided we can show that
\begin{equation}
\label{eqn: remaining vanishing}
\chi_{\infty, \infty}^{\underline{x}, \underline{y}}(R^{A_\infty}_\perfd \otimes^\myL_SG) =0.
\end{equation}

Note that  $\lambda_{\infty,\infty}^{\underline{x},\underline{y}}\big(H_i(R^{A_\infty}_\perfd \otimes^\myL_SG)\big)<\infty$ for all $i$ by examining the long exact sequence of \autoref{eqn: isogeny} derived tensor with $R^{A_\infty}_\perfd$ (thus in particular, $\chi_{\infty, \infty}^{\underline{x}, \underline{y}}(R^{A_\infty}_\perfd \otimes^\myL_SG)$ is well-defined). We will prove \autoref{eqn: remaining vanishing} via ``reduction to the diagonal".
\end{discussion}

For the ease of presentation, we set up some notations.

\begin{setting}
\label{setting: step 2}
We continue to use the notations as in \autoref{setting: step 1} and \autoref{setting: normalized length}. We write $S=T/(p-f)$ where $T= W(k)[[z_1,\dots,z_n]]$ and $f\in (p,z_1,\dots,z_n)^2$. For two $T$-modules $M$ and $N$, we will use $M\widehat{\otimes}_{W(k)}^\myL N$ to denote the derived completion of $M\otimes_{W(k)}^\myL N$ with respect to the ideal $(p, z_i\otimes 1, 1\otimes z_j \mid 1\leq i, j \leq n) \subseteq T\otimes_{W(k)}T$. In particular, $T\widehat{\otimes}_{W(k)}^\myL T$ is isomorphic to a power series ring over $W(k)$ in $2n$ variables, which can be identified with $z_i\otimes 1$ and $1\otimes z_j$ for $1\leq i, j \leq n$. We write $\delta_i := z_i\otimes 1 - 1\otimes z_i \in T\widehat{\otimes}_{W(k)}^\myL T$. Fix $h\notin Q$ (this will be a power of the discriminant of $B\to S/Q$ as in \autoref{discuss: strategy of proof}) and fix a map $T\to A_{\inf}(R^{A_\infty}_\perfd)$ as in \autoref{lem: Tor-independent}, we consider
$$A_{\inf}(R^{A_\infty}_\perfd) \widehat{\otimes}^\myL_{W(k)} S/(Q+h)$$
as a $T\widehat{\otimes}^\myL_{W(k)} T$-algebra. Note that $(T\widehat{\otimes}^\myL_{W(k)} T)/^\myL (\delta_1,\dots,\delta_n) \cong T$ and by \autoref{lem: Tor-independent},
\begin{align}
\label{eqn: eqn in the last setup}
& \big(A_{\inf}(R^{A_\infty}_\perfd) \widehat{\otimes}^\myL_{W(k)} S/(Q+h)\big) /^{\myL}(\delta_1,\dots,\delta_n) \\
 = & \Kos\big(\delta_1,\dots,\delta_n; \big(A_{\inf}(R^{A_\infty}_\perfd) \widehat{\otimes}^\myL_{W(k)} S/(Q+h)\big) \big) \notag \\
 \cong & A_{\inf}(R^{A_\infty}_\perfd) {\otimes}^\myL_{T} S/(Q+h) \notag \\ 
 =  & R^{A_\infty}_\perfd \otimes^\myL_S S/(Q+h).    \notag
\end{align}
\end{setting}

\begin{remark}
We note that $A_{\inf}(R^{A_\infty}_\perfd) \widehat{\otimes}^\myL_{W(k)}S/(Q+h)$ is a discrete (local) ring. In fact, for any finitely generated $T$-module $N$, $A_{\inf}(R^{A_\infty}_\perfd) \widehat{\otimes}^\myL_{W(k)} N$ is discrete. To see this, we may assume that $N$ is either $p$-torsion free or $p$-torsion by filtering $N$. In either case, it suffices to show that $(R^{A_\infty}_\perfd)^\flat\widehat{\otimes}_k N'$ is discrete where $N'$ is a finitely generated $T/p$-module. By Noetherianness of $N'$, $z_1,\dots,z_n$ is a weakly pro-regular sequence on $N'$; and since $(R^{A_\infty}_\perfd)^\flat$ is perfect, the images of $z_1,\dots,z_n$ form a pro-regular sequence (and hence a weakly pro-regular sequence) on $(R^{A_\infty}_\perfd)^\flat$.\footnote{If $R$ is perfect and $f_1,\dots,f_n$ is any sequence in $R$, then $(f_1^{pm},\dots,f_i^{pm}):f_{i+1}^{pm}\subseteq (f_1^m,\dots,f_i^m):f_{i+1}^{pm-m}$.} By examining the Koszul complex on $z^m_1\otimes1,\dots, z^m_n\otimes 1, 1\otimes z^m_1,\dots,1\otimes z^m_n$ it is then easy to see that $z_1\otimes1,\dots, z_n\otimes 1, 1\otimes z_1,\dots,1\otimes z_n$ is a weakly pro-regular sequence on $(R^{A_\infty}_\perfd)^\flat{\otimes}_k N'$. Thus its derived completion is discrete. (See \cite{SchenzelProregularSequences2003} for various facts we used about (weakly) pro-regular sequences and derived completion.)
\end{remark}

Our idea for proving \autoref{eqn: remaining vanishing} is to apply \autoref{eqn: eqn in the last setup} to reduce the computation of $\chi_{\infty,\infty}^{\underline{x},\underline{y}}$ to a Koszul complex, and then show that under \autoref{setting: step 2}, one of the generators of $(\delta_1,\dots,\delta_n)$ is redundant (up to a power) in $A_{\inf}(R^{A_\infty}_\perfd) \widehat{\otimes}^\myL_{W(k)} S/(Q+h)$. In order to establish this, we need some preparations on valuation theory (which allow us to prove a stronger result).

For a valuation $v$ on a ring $R$, we write $\Gamma_v$ and $k_v$ for the value group and the residue field of the induced valuation on the fraction field of $R/Q$, where $Q:= v^{-1}(\infty)$ is the support of $v$ on $R$. We recall two basic observations on the behavior of $\Gamma_v$ and $k_v$: %when we extend the valuation:
\begin{itemize}
    \item Suppose $R\to R'$ is an integral extension and $w$ is a valuation on $R'$ that extends a valuation $v$ on $R$. Then $\Gamma_w/\Gamma_v$ is torsion and $k_w/k_v$ is an algebraic extension.
    \item Suppose $v$ is a rank one valuation on $R$ and $J$ is a finitely generated ideal of $R$ that is contained in the center of $v$ (i.e., $v(J)>0$). Let $w$ be the valuation defined on $R^{\wedge J}$ by $w(\lim_i r_i):= \lim_i v(r_i)$. Then $\Gamma_w=\Gamma_v$ and $k_w=k_v$. 
\end{itemize}

We record two results on the behavior of value groups and residue fields under ring operations (these results should be well-known to experts). First, we study them under perfectoidizations.

\begin{lemma}
\label{lem: extend valuation perfectoidization}
Let $R$ be a ring such that $R_\perfd$ is discrete (and hence a perfectoid ring). Suppose $w$ is a rank one valuation on $R_\perfd$ so that $w(p)> 0$ and let $v$ be its restriction to $R$. Then $\Gamma_w/\Gamma_v$ is torsion and $k_w/k_v$ is an algebraic extension.  
\end{lemma}
\begin{proof}
If $w(p)=\infty$, then $w$ is a valuation on $R_\perfd/p$ and the conclusion follows from that $R/p\to R_\perfd/p$ is integral. Henceforth we assume $w(p)<\infty$. Let $Q\in\Spec(R)$ be the support of $v$ and let $F$ be the fraction field of $R/Q$. Similarly, let $Q'\in\Spec(R_\perfd)$ be the support of $w$ and let $F'$ be the fraction field of $R_\perfd/Q'$. Fix a $p$-completed algebraic closure $\widehat{\overline{F'}}$ of $F'$ and let $\widehat{\overline{F}}$ be the $p$-adic closure of the algebraic closure of $F$ inside $\widehat{\overline{F'}}$. Then $\sO_{\widehat{\overline{F}}}$ is a rank one perfectoid valuation ring and thus we have 
$$R\to R_\perfd\to \sO_{\widehat{\overline{F}}}\to \sO_{\widehat{\overline{F'}}}$$ by the universal property of perfectoidization. In particular, the map $R_\perfd\to \sO_{\widehat{\overline{F}}}$ factors through $R_\perfd/Q'$ and thus we have 
$$F\to F'\to \widehat{\overline{F}} \to \widehat{\overline{F'}}$$
as maps of valued fields. In particular, if $\overline{w}$ is the rank one valuation on $\widehat{\overline{F}}$ with valuation ring $\sO_{\widehat{\overline{F}}}$, then $w$ and $v$ are its restrictions to $F'$ and $F$ respectively. Since $\Gamma_{\overline{w}}/\Gamma_v$ is torsion and $k_{\overline{w}}/k_v$ is algebraic, it follows that $\Gamma_{w}/\Gamma_v$ is torsion and $k_{w}/k_v$ is algebraic. 
\end{proof}

The next lemma is well-known to experts and follows from more general results in \cite[Corollary 8.19]{KedlayaWittVectors2013}. We include a short proof in our context.

\begin{lemma}
\label{lem: valuation on W(R)}
Let $R$ be a perfect ring of characteristic $p$ and let $w$ be a rank one valuation on $W(R)$ such that $w(p)>0$. Then $v(r):=w([r])$ defines a rank at most one valuation on $R$ with natural inclusions $\Gamma_v\subseteq \Gamma_w$ and $k_v\subseteq k_w$. Moreover, 
$$\dim_\Q(\Gamma_w/\Gamma_v \otimes_\Z\Q) + \trdeg_{k_v}k_w\leq 1.$$
\end{lemma}
\begin{proof}
%First note that if $p$ is in the support of $v$, then $w=v$ is a valuation on $W(R)/pW(R)=R$ and the conclusion is obvious. 
Note that $R':=W(R)[p^{1/p^\infty}]^{\wedge p}$ is a perfectoid ring and $W(R)\to R'$ is a $p$-completed integral extension. Thus after extending $w$ to $w'$ on $R'$ (by setting $w'(p^{1/p^e})=w(p)/p^e$), we know that $\Gamma_{w'}/\Gamma_{w}$ is torsion and $k_{w'}/k_w$ is algebraic. 

By the tilting correspondence \cite[Lemma 3.4 (iii)]{ScholzePerfectoidSpaces2012}, $w'$ induces a rank one valuation $w''$ on $R'^{\flat}\cong R[t^{1/p^\infty}]^{\wedge t}$ so that $\Gamma_{w'}=\Gamma_{w''}$ and $k_{w'}=k_{w''}$. It is straightforward to check that $w''|_R=v$ where $v$ is defined as in the lemma (thus in particular, $v$ is a rank at most one\footnote{It is possible that $v$ is the trivial valuation (e.g., $R=\mathbb F_p$ and $w$ is the $p$-adic valuation on $\mathbb Z_p$).} valuation on $R$). Let $v':=w''|_{R[t]}$. Then as $R[t]\to R[t^{1/p^\infty}]^{\wedge t}$ is a $t$-completed integral extension and that $v'(t)=w''(t)=w'(p)>0$ (so $t$ is in the center), we know that $\Gamma_{w''}/\Gamma_{v'}$ is torsion and ${k_{v'}}\to k_{w''}$ is algebraic. Finally, by the relative version of the Abhyankar inequality \cite[VI.10.3,
Corollary 1]{Bourbaki1998}, we have 
$$\dim_\Q(\Gamma_{v'}/\Gamma_v \otimes_\Z\Q) + \trdeg_{k_v}k_{v'}\leq 1.$$
Thus the same holds upon replacing $\Gamma_{v'}$ by $\Gamma_w$ and $k_{v'}$ by $k_w$. %This completes the proof of the lemma.
\end{proof}

The next result is a key ingredient in our proof of \autoref{eqn: remaining vanishing}.

\begin{proposition}
\label{prop: analytic spread at most n-1}
With notations as in \autoref{setting: step 2}. Suppose $\widetilde{\delta}_1,\dots, \widetilde{\delta}_n$ are general $W(k)$-linear combinations of $\delta_1,\dots,\delta_n$. Then the ideal $(\delta_1,\dots,\delta_n)\big(A_{\inf}(R^{A_\infty}_\perfd) \widehat{\otimes}^\myL_{W(k)} S/(Q+h)\big)$ is integral over $(\widetilde{\delta}_1,\dots, \widetilde{\delta}_{n-1})\big(A_{\inf}(R^{A_\infty}_\perfd) \widehat{\otimes}^\myL_{W(k)}S/(Q+h)\big)$. In particular, there exists an~integer $L>0$ so that 
$$\widetilde{\delta}_n^L \in (\widetilde{\delta}_1,\dots, \widetilde{\delta}_{n-1})\big(A_{\inf}(R^{A_\infty}_\perfd) \widehat{\otimes}^\myL_{W(k)}S/(Q+h)\big).$$
\end{proposition}
\begin{proof}
We use $(U,\m_U, k)$ to denote the local ring $A_{\inf}(R^{A_\infty}_\perfd) \widehat{\otimes}^\myL_{W(k)}S/(Q+h)$ in what follows. We will examine the rank one valuations on $U$.

\begin{claim}
\label{clm: clm in analytic spread prop}
For any rank one
valuation $v$ on $U$ centered at $\m_U$, we have $\trdeg_kk_v<n-1$.
\end{claim}
\begin{proof}[Proof of Claim]
Suppose $\trdeg_kk_v\geq n-1$. Then we have 
$$\dim_\Q(\Gamma_v\otimes_\Z\Q) + \trdeg_kk_v \geq n.$$

Let $v_0$ be the restriction of $v$ to $A_{\inf}(R^{A_\infty}_\perfd) {\otimes}^\myL_{W(k)}S/(Q+h)$ and let $Q'$ be the support of (the restriction of) $v_0$ on $S/(Q+h)$. We thus obtain a valuation $v_0$ on $A_{\inf}(R^{A_\infty}_\perfd) {\otimes}^\myL_{W(k)} S/Q'$ with the same value group and residue field as $v$, in particular
$$\dim_\Q(\Gamma_{v_0}\otimes_\Z\Q) + \trdeg_kk_{v_0} \geq n.$$

Note that $\dim(S/Q')\leq \dim(S/(Q+h))=c-1$. If $p\notin Q'$, then let $W(k)[[t_2,\dots,t_{c-1}]] \to S/Q'$ be a module-finite map. If $p\in Q'$, then let $k[[t_1,\dots,t_{c-1}]]\to S/Q'$ be a module-finite map. In either case, let $v_1$ be the restriction of $v_0$ to
$$A_{\inf}(R^{A_\infty}_\perfd) {\otimes}^\myL_{W(k)} W(k)[[t_2,\dots,t_{c-1}]] \,\ \text{ or } \,\ A_{\inf}(R^{A_\infty}_\perfd) {\otimes}^\myL_{W(k)} k[[t_1,\dots,t_{c-1}]]$$
respectively. In both cases we have 
$$\dim_\Q(\Gamma_{v_1}\otimes_\Z\Q) + \trdeg_kk_{v_1} \geq n.$$

Further let $v_2$ be the restriction of $v_1$ to
$$A_{\inf}(R^{A_\infty}_\perfd) {\otimes}^\myL_{W(k)} W(k)[t_2,\dots,t_{c-1}] \,\ \text{ or } \,\ A_{\inf}(R^{A_\infty}_\perfd) {\otimes}^\myL_{W(k)} k[t_1,\dots,t_{c-1}].$$
Then $v_2$ and $v_1$ have the same value group and residue field and so
$$\dim_\Q(\Gamma_{v_2}\otimes_\Z\Q) + \trdeg_kk_{v_2} \geq n.$$
By the relative version of the Abhyankar inequality (see \cite[VI.10.3,
Corollary 1]{Bourbaki1998}) and \autoref{lem: valuation on W(R)}, we then obtain that in either case, there is a rank one valuation $v'$ on $(R^{A_\infty}_\perfd)^\flat$ centered at its maximal ideal so that 
$$\dim_\Q(\Gamma_{v'}\otimes_\Z\Q) + \trdeg_kk_{v'} \geq n-(c-1).$$
It follows from the tilting correspondence (\cite[Lemma 3.4 (iii)]{ScholzePerfectoidSpaces2012}) that there is a rank one valuation $w$ on $R^{A_\infty}_\perfd$ centered at its maximal ideal so that 
$$\dim_\Q(\Gamma_{w}\otimes_\Z\Q) + \trdeg_kk_{w} \geq n-(c-1).$$

Finally, let $w'$ be the restriction of $w$ to $R$, since $R\to R\otimes_AA_\infty$ is a $p$-completed integral extension and $(R\otimes_AA_\infty)_\perfd=R^{A_\infty}_\perfd$ is discrete, by \autoref{lem: extend valuation perfectoidization} we obtain that 
$$\dim_\Q(\Gamma_{w'}\otimes_\Z\Q) + \trdeg_kk_{w'} \geq n-(c-1).$$
But this contradicts the classical Abhyankar's inequality \cite[Theorem 1]{AbhyankarValuations1956}, namely that $\dim_\Q(\Gamma_{w'}\otimes_\Z\Q) + \trdeg_kk_{w'} \leq \dim(R)=d =n-c.$
\end{proof}

We now come back to the proof of the proposition. Let $J:=(\delta_1,\dots,\delta_n)U$. By \autoref{clm: clm in analytic spread prop} and \autoref{lem: analytic spread n} below, we know that the finitely generated $k$-algebra 
$$\mathcal{F}_J(U):=k \oplus J/\m_UJ \oplus J^2/\m_UJ^2\oplus \cdots$$
has dimension at most $n-1$. Since $k$ is infinite, it follows that the images of $\widetilde{\delta}_1,\dots, \widetilde{\delta}_{n-1}$ in $\mathcal{F}_J(U)$ form a homogeneous system of parameters. Hence for $L\gg0$, 
$$J^L \subseteq (\widetilde{\delta}_1,\dots, \widetilde{\delta}_{n-1})J^{L-1} + \m_U J^L$$
and thus $J^L=(\widetilde{\delta}_1,\dots, \widetilde{\delta}_{n-1})J^{L-1}$ by Nakayama's lemma. This proves that $J$ is integral over $(\widetilde{\delta}_1,\dots, \widetilde{\delta}_{n-1})U$ and the last assertion is obvious from this.
\end{proof}

The following simple lemma is used above. 

\begin{lemma}
\label{lem: analytic spread n}
Let $(R,\m_R, k)$ be a (not necessarily Noetherian) local ring and let $J=(f_1,\dots,f_n)\subseteq R$ be a finitely generated ideal so that $\sqrt{J}=\m_R$. Suppose the natural map 
$$R/\m_R[T_1,\dots,T_n] \twoheadrightarrow \mathcal{F}_J(R):= R/\m_R \oplus J/\m_RJ \oplus J^2/\m_RJ^2 \oplus \cdots $$
sending $T_i$ to $\overline{f}_i\in J/\m_R J$ is an isomorphism. Then there is a rank one valuation $v$ of $R$ centered at $\m_R$ so that 
$\trdeg_kk_v\geq n-1$.
\end{lemma}
\begin{proof}
Consider the commutative diagram: 
\[\xymatrix{
\Proj \mathcal{F}_J(R) \cong \mathbb{P}_{k}^{n-1} \ar@{^{(}->}[r] \ar[d] & Y := \Bl_J(R) \ar[d] \\
\Spec(k) \ar@{^{(}->}[r] &  \Spec(R).
}\]
Let $\eta$ be the generic point of the irreducible subscheme $\mathbb{P}_{k}^{n-1}$ inside $Y$. Then we have a map of local rings:
$$(R,\m_R) \to (\sO_{Y,\eta},\m_{\eta}).$$
Let $v$ be any valuation of $\sO_{Y,\eta}$ centered at $\m_{\eta}$ and supported at some prime ideal $\mathfrak{Q}\neq \m_\eta$ of $\sO_{Y,\eta}$. Let $Q$ be the contraction of $\mathfrak{Q}$ to $R$. Then $v$ is a valuation of $R$ centered at $\m_R$ and supported at $Q$ (note that the fraction field of $R/Q$ agrees with the fraction field of $\sO_{Y,\eta}/\mathfrak{Q}$). It is clear that $\trdeg_kk_v\geq n-1$ since by construction $k_v$ is an extension of the function field of $\mathbb{P}_k^{n-1}$. Also note that by construction, $v(f_i)=v(f_j)$ for all $i$ and $j$ since $T_i/T_j\in \sO_{Y,\eta}$ is a unit. Thus $f_i=u_i\cdot f \in V$ where $V$ is the valuation ring of $v$ and $u_i\in V^\times$. Note that 
$$\mathfrak{p}:= \cap_n(fV)^n\subseteq \sqrt{fV}=: \mathfrak{q}$$
are two intermediate prime ideals of $V$. Replacing $v$ by the rank one valuation corresponding to the valuation ring $(V/\mathfrak{p})_\mathfrak{q}$, it is easy to see that $v$ is still a valuation on $R$ centered at $\m_R$ with $\trdeg_kk_v\geq n-1$. 
\end{proof}

Finally, we have all the ingredients to prove the main result that implies Serre positivity.

\begin{theorem}
\label{thm: vanishing and completion of proof}
With notations as in \autoref{discuss: strategy of proof}. %Suppose $G$ is a module over $S/(Q+h) \otimes_BB_\infty$ so that $\lambda_{\infty,\infty}^{\underline{x},\underline{y}}\big(H_i(R^{A_\infty}_\perfd \otimes^\myL_SG)\big)<\infty$ for all $i$. 
We have
$$\chi_{\infty, \infty}^{\underline{x}, \underline{y}}(R^{A_\infty}_\perfd \otimes^\myL_SG) =0.$$
In particular, we have 
$$\chi(S/P \otimes_S^\myL S/Q)>0.$$
\end{theorem}
\begin{proof}
The ``In particular" part of the theorem follows from \autoref{discuss: strategy of proof}, provided we can show the first statement of the theorem. 

To prove the first statement, we follow the notations as in \autoref{setting: step 2}. Let $\widetilde{\delta}_1,\dots,\widetilde{\delta}_n$ be general $W(k)$-linear combinations of $\delta_1,\dots,\delta_n$ that satisfy the conclusions of \autoref{prop: analytic spread at most n-1}, and we fix an integer $L>0$ so that 
\begin{equation}
\label{eqn: gen up to radical n-1 elts}
\widetilde{\delta}_n^L \in (\widetilde{\delta}_1,\dots,\widetilde{\delta}_{n-1})\big(A_{\inf}(R^{A_\infty}_\perfd) \widehat{\otimes}^\myL_{W(k)} S/(Q+h)\big)
\end{equation}
by \autoref{prop: analytic spread at most n-1}.  Now we have
\begin{align*}
R^{A_\infty}_\perfd \otimes^\myL_SG & \cong A_{\inf}(R^{A_\infty}_\perfd) \otimes^\myL_TG\\
& \cong \Kos(\delta_1,\dots,\delta_n; A_{\inf}(R^{A_\infty}_\perfd) \widehat{\otimes}^\myL_{W(k)}G) \\
& \cong  \Kos(\widetilde{\delta}_1,\dots,\widetilde{\delta}_n; A_{\inf}(R^{A_\infty}_\perfd) \widehat{\otimes}^\myL_{W(k)}G
) 
\end{align*}
where the first equality follows from \autoref{lem: Tor-independent}. In particular, we have 
\begin{align}
\label{eqn: reduction to diagonal}
\chi_{\infty, \infty}^{\underline{x}, \underline{y}}(R^{A_\infty}_\perfd \otimes^\myL_SG) & = \chi_{\infty, \infty}^{\underline{x}, \underline{y}}\big(\Kos(\widetilde{\delta}_1,\dots,\widetilde{\delta}_n; A_{\inf}(R^{A_\infty}_\perfd) \widehat{\otimes}^\myL_{W(k)}G
)\big) \notag \\
& = \frac{1}{L} \cdot \chi_{\infty, \infty}^{\underline{x}, \underline{y}}\big(\Kos(\widetilde{\delta}_1,\dots,\widetilde{\delta}_{n-1},\widetilde{\delta}_n^L; A_{\inf}(R^{A_\infty}_\perfd) \widehat{\otimes}^\myL_{W(k)}G
)\big)
\end{align}
where the second equality follows by filtering the Koszul complex. Note that by filtering the Koszul complex and examining the long exact sequence of Koszul homologies, we also know that
$$\lambda_{\infty, \infty}^{\underline{x}, \underline{y}}\big(H_i\big(\Kos(\widetilde{\delta}_1,\dots,\widetilde{\delta}_{n-1}, \widetilde{\delta}_n^L; A_{\inf}(R^{A_\infty}_\perfd) \widehat{\otimes}^\myL_{W(k)}G)\big)\big) <\infty$$
for all $i$ (since $\lambda_{\infty,\infty}^{\underline{x},\underline{y}}\big(H_i(R^{A_\infty}_\perfd \otimes^\myL_SG)\big)<\infty$ for all $i$, see the last paragraph of \autoref{discuss: strategy of proof}). Now since $G$ is a module over $S/(Q+h)$, by \autoref{eqn: gen up to radical n-1 elts} we know that 
\begin{align*}
& \,\ H_i\big(\Kos(\widetilde{\delta}_1,\dots,\widetilde{\delta}_{n-1},\widetilde{\delta}_n^L; A_{\inf}(R^{A_\infty}_\perfd) \widehat{\otimes}^\myL_{W(k)}G)\big) \\
\cong & \,\ H_i\big(\Kos(\widetilde{\delta}_1,\dots,\widetilde{\delta}_{n-1}, 0; A_{\inf}(R^{A_\infty}_\perfd) \widehat{\otimes}^\myL_{W(k)}G)\big) \\
\cong & \,\ H_i\big(\Kos(\widetilde{\delta}_1,\dots,\widetilde{\delta}_{n-1}; A_{\inf}(R^{A_\infty}_\perfd) \widehat{\otimes}^\myL_{W(k)}G)\big)  \oplus H_{i-1}\big(\Kos(\widetilde{\delta}_1,\dots,\widetilde{\delta}_{n-1}; A_{\inf}(R^{A_\infty}_\perfd) \widehat{\otimes}^\myL_{W(k)}G)\big)
\end{align*}
for all $i$. In particular, this implies that 
$$\lambda_{\infty, \infty}^{\underline{x}, \underline{y}}\big(H_i\big(\Kos(\widetilde{\delta}_1,\dots,\widetilde{\delta}_{n-1}; A_{\inf}(R^{A_\infty}_\perfd) \widehat{\otimes}^\myL_{W(k)}G)\big)\big) <\infty$$
for all $i$, and after taking the alternating sum of the normalized length we have 
\begin{align}
\label{eqn: Koszul cancelation}
& \,\ \chi_{\infty, \infty}^{\underline{x}, \underline{y}}\big(\Kos(\widetilde{\delta}_1,\dots,\widetilde{\delta}_{n-1},\widetilde{\delta}_n^L; A_{\inf}(R^{A_\infty}_\perfd) \widehat{\otimes}^\myL_{W(k)}G
)\big) \\
= & \,\  \sum_i(-1)^i\lambda_{\infty, \infty}^{\underline{x}, \underline{y}}\big(H_i\big(\Kos(\widetilde{\delta}_1,\dots,\widetilde{\delta}_{n-1}, \widetilde{\delta}_n^L; A_{\inf}(R^{A_\infty}_\perfd) \widehat{\otimes}^\myL_{W(k)}G)\big)\big) \notag\\
= & \,\ \sum_i(-1)^i\lambda_{\infty, \infty}^{\underline{x}, \underline{y}}\big(H_i\big(\Kos(\widetilde{\delta}_1,\dots,\widetilde{\delta}_{n-1}; A_{\inf}(R^{A_\infty}_\perfd) \widehat{\otimes}^\myL_{W(k)}G)\big)\big) \notag \\
&  \,\ \,\ + \sum_i(-1)^i\lambda_{\infty, \infty}^{\underline{x}, \underline{y}}\big(H_{i-1}\big(\Kos(\widetilde{\delta}_1,\dots,\widetilde{\delta}_{n-1}; A_{\inf}(R^{A_\infty}_\perfd) \widehat{\otimes}^\myL_{W(k)}G)\big)\big)  \notag\\
= & \,\ 0. \notag
\end{align}
Putting \autoref{eqn: reduction to diagonal} and \autoref{eqn: Koszul cancelation} together completes the proof of the theorem.
\end{proof}

\begin{proof}[Proof of \autoref{theorem: Serre positivity}]
Without loss of generality, we may assume $S$ is complete with an algebraically closed residue field. If both $S/P$ and $S/Q$ have characteristic $p$, then the result was known \cite{MalliavinPositivity1972} (alternatively, both $S/P$ and $S/Q$ admit lim Cohen--Macaulay sequences by \cite[Theorem 5.4]{BhattHochsterMaLimCM2026} so the result also follows from \cite[Theorem 6.1]{BhattHochsterMaLimCM2026}). If $S/P$ is $p$-torsion free and $S/Q$ has characteristic $p$, then the result follows from \autoref{thm: S/Q has char p}. Finally, if both $S/P$ and $S/Q$ are $p$-torsion free, then the result follows from \autoref{thm: vanishing and completion of proof}.
\end{proof}

\bibliographystyle{skalpha}
\bibliography{MainBib}
\end{document}